\documentclass[journal,10pt]{IEEEtran}

\usepackage{graphicx}  

\usepackage{subfigure} 

\usepackage{amsmath}   

\usepackage{amssymb}

\usepackage{rotating}
\usepackage{cite}
\usepackage{rotating}
\usepackage{pstricks}

\newtheorem{theorem}{Theorem}

\def\x{{\mathbf x}}
\def\y{{\mathbf y}}
\def\u{{\mathbf u}}

\def\C{{\mathbf C}}
\def\D{{\mathbf D}}

\def\A{{\mathbf A}}
\def\B{{\mathbf B}}

\def\N{{\mathbf N}}
\def\F{{\mathbf F}}
\def\Q{{\mathbf Q}}
\def\K{{\mathbf K}}
\def\R{{\mathbf R}}
\def\S{{\mathbf S}}
\def\w{{\mathbf w}}

\def\n{{\mathbf n}}
\def\H{{\mathbf H}}

\def\I{{\mathbf I}}
\def\Z{{\mathbf Z}}

\def\s{{\mathbf s}}
\def\Q{{\mathbf Q}}
\def\Z{{\mathbf Z}}

\def\e{{\mathbf e}}
\def\vec{{\text{vec}}}

\begin{document}

\title{MMSE Estimation Under \\ Gaussian Mixture Statistics}
\title{The Linear Model under Mixed Gaussian Inputs:\\Designing the Transfer Matrix}

\author{John T. Fl{\aa}m, Dave Zachariah, Mikko Vehkaper\"{a} and Saikat Chatterjee
\thanks{John T. Fl{\aa}m is with the Department of Electronics and Telecommunications, NTNU-Norwegian University of Science and Technology, Trondheim, Norway. Email: flam@iet.ntnu.no.

 Dave Zachariah, Mikko Vehkaper\"{a} and Saikat Chatterjee are with the School of Electrical Engineering, KTH-Royal Institute of Technology, Sweden. Emails: davez@kth.se, mikkov@kth.se, sach@kth.se.
}}

\maketitle

\begin{abstract}

Suppose a linear model $\mathbf{y}=\mathbf{H}\mathbf{x}+\mathbf{n}$, where inputs $\mathbf{x,n}$ are independent Gaussian mixtures. 
The problem is to design the transfer matrix $\mathbf{\H}$ so as to minimize the mean square error (MSE) when estimating $\mathbf{x}$ from $\mathbf{y}$.  
This problem has important applications, but faces at least three hurdles. Firstly, even for a fixed $\H$, the minimum MSE (MMSE) has no analytical form. Secondly, the MMSE is generally not convex in $\H$. Thirdly, derivatives of the MMSE w.r.t. $\H$ are hard to obtain. 
This paper casts the problem as a stochastic program and invokes gradient methods. 

The study is motivated by two applications in signal processing. One
concerns the choice of error-reducing precoders; the other deals with
selection of pilot matrices for channel estimation. In either setting, our
numerical results indicate improved estimation accuracy - markedly better
than those obtained by optimal design based on standard linear estimators.

Some implications of the non-convexities of the MMSE are noteworthy, yet, to our knowledge,
not well known. For example, there are cases in which more pilot power is
detrimental for channel estimation. This paper explains why. 
%


\end{abstract}

%
\IEEEpeerreviewmaketitle

\section{Problem statement}\label{intro}

Consider the following linear system 
\begin{equation}
\mathbf{y=Hx+n}.
\label{linmod}
\end{equation}
Here $\mathbf{y}$ is a vector of observations, and $\mathbf{x}$ and $\mathbf{n}$ are mutually independent random vectors with known Gaussian Mixture (GM) distributions:
\begin{align}
&\mathbf{x}\sim \hspace{-0.1cm}\sum_{k\in\mathcal{K}}p_k \mathcal{N}\hspace{-0.1cm}\left(\mathbf{u}^{(k)}_{\mathbf{x}},\mathbf{C}^{(k)}_{\mathbf{x}\mathbf{x}}\right)\label{x_mix}\\
 &\mathbf{n}\sim \hspace{-0.1cm}\sum_{l\in\mathcal{L}}q_l \mathcal{N}\hspace{-0.1cm}\left(\mathbf{u}^{(l)}_{\mathbf{n}},\mathbf{C}^{(l)}_{\mathbf{n}\mathbf{n}}\right).
\label{n_mix}
\end{align}
In this work, we assume that $\H$ is a transfer matrix that we are at liberty to \textit{design}, typically under some constraints. Specifically, our objective is to design $\H$ such that $\x$ can be estimated from $\y$ with minimum mean square error (MMSE).
The MMSE, for a \textit{fixed} $\H$, is by definition \cite{151045}
\begin{align}
\text{MMSE}&\triangleq E \left\{ \| \mathbf{x} -\u_{\x|\y} \|_2^{2}\right\}\nonumber\\
&=\iint \| \mathbf{x} -\u_{\x|\y} \|_{2}^{2}f(\x,\y)d\x d\y\label{defmmse}.
\end{align}
Here, $\left\|\cdot\right\|_2$ denotes the 2-norm, $f(\x,\y)$ is the joint probability density function (PDF) of $(\x,\y)$,
\begin{align}
\u_{\x|\y}\triangleq E\left\{\x|\y\right\}=\int \x f(\x|\y)d\x \label{defmmseest}
\end{align}
is the MMSE estimator, and $f(\x|\y)$ is the PDF of $\x$ given $\y$. The MMSE in equation \eqref{defmmse} depends on $\H$ both through $\u_{\x|\y}$  and $f(\x,\y)$. Our objective is to solve the following optimization problem
\begin{align}
&\min_{\H \in \mathbb{H} } \hspace{0.2cm}\text{MMSE}\label{theproblem}, 
\end{align}
where $\mathbb{H}$ denotes a prescribed set of matrices that $\H$ must belong to.
Solving this optimization problem is not straightforward. In particular, three hurdles stand out. Firstly, with (\ref{x_mix}) and \eqref{n_mix} as inputs to (\ref{linmod}), the MMSE in (\ref{defmmse}) has no analytical closed form \cite{flamchatterjeekansanen}. 
Thus, the effect of \textit{any} matrix $\mathbf{H}$, in terms of MMSE, cannot be evaluated exactly. Secondly, the MMSE is not convex in $\H$. Thirdly, the first and second order derivatives of the MMSE w.r.t $\H$ cannot be calculated exactly, and accurate approximations are hard to obtain. For these reasons, and in order to make progress, we cast the problem as a stochastic program and invoke the Robbins-Monro algorithm\cite{robbins1951, kushner2003stochastic}. 
Very briefly our approach goes as follows: We draw samples from $\x$ and $\n$ and use these to compute
stochastic gradients of the MMSE. These feed into an
iterative gradient method that involves projection.


The contributions of the paper are several:
\begin{itemize}
	\item As always, for greater accuracy, its preferable to use gradients
instead of finite difference approximations. For this reason the paper spells out a
formula for exact realization of stochastic gradients. Accordingly, the
Robbins-Monro algorithm comes to replace the Kiefer-Wolfowitz procedure.
\item In the design phase, we exploit the known input statistics and update $\H$ based on samples of the inputs $(\x,\n)$, instead of output $\y$. This yields a closed form stochastic gradient, and we prove that it is unbiased. 
\item Numerical experiments indicate that our method has far better accuracy than methods which
proceed via linear estimators. The main reason is that the optimal estimator, used here, is non-linear.
\item It turns out that the non-convexities of the MMSE may have practical implications
that deserve being better known. Specifically, in channel estimation, it can
be harmful to increase the power of the pilot signal. This paper offers an explanation.
\end{itemize}

Clearly, in many practical problems, the quantities in (\ref{linmod}) are complex-valued. Throughout this paper, however, they will all be assumed real.
For the analysis, this assumption introduces no loss of generality, as the real and imaginary parts of (\ref{linmod}) can always be treated separately.

The paper is organized as follows. The next section outlines two applications. It
also specifies the Gaussian mixtures and motivates their use. Section \ref{A simple example}
illustrates the problem by means of a simple example. Section \ref{An equivalent optimization problem} spells out
problem (\ref{theproblem}) in full detail. Section \ref{Kiefer-Wolfowitz precoding} reviews how the
Robbins-Monro method applies. Numerical results are provided in Section \ref{Numerical results}.
Section \ref{concl} concludes. A large part of the detailed analysis, concerning
stochastic gradients, is deferred to the appendix.

\section{Background and Motivation}\label{Background and Motivation}

The above described matrix design problem appears in various applications of interest. Next we present two of these, which are of particular interest to the signal processing community. Then we will explain and motivate the GM input statistics.

\subsection {Linear precoder design}\label{Linear precoding}
Consider a linear system model
\begin{equation}
\y=\underbrace{\B\F}_{\H}\x+\n,
\label{linmod5}
\end{equation}
where $\B$ is a known matrix, and $\F$ is a precoder matrix to be designed such that the mean square error (MSE) when estimating $\mathbf{x}$ from $\mathbf{y}$ becomes as small as possible. The vector $\n$ is random noise. If $\x$ and $\n$ are independent and GM distributed, this is a matrix design problem as described in Section \ref{intro}, where $\F$ is the design parameter. 
A typical constraint is to require that $\F\x$ cannot exceed a certain average power, i.e. $E\left\|\F\x\right\|_2^2\leq \gamma$. Together with the nature of $\B$, this determines $\mathbb{H}$ in \eqref{theproblem}.

A linear model with known transfer matrix $\H$ and GM distributed inputs is frequently assumed within speech and image processing. In these applications, the signal of interest often exhibits multi modal behavior. That feature can be reasonably explained by assuming an underlying GM distribution. The noise is often modeled as Gaussian, which is a special case of a GM. Conveniently, with (\ref{x_mix}) and \eqref{n_mix} as inputs, the MMSE estimator in (\ref{defmmseest}) has a closed analytical form for any given $\H$. Selected works exploiting this include \cite{4808405,kundusaikat,kundusaikatsr,1597257,6104390}. However, none of these works study MSE reducing precoders. These have the potential to significantly improve the estimation accuracy, and should therefore be of interest. 

\subsection {Pilot signal design}\label{MIMO channel estimation}
Consider a multiple-input-multiple-output (MIMO) communication model
\begin{equation}
\mathbf{z=As+n},
\label{linmod2}
\end{equation}
where $\A$ is a random channel matrix that we wish to estimate with as small MSE as possible, and $\s$ is a pilot signal to be designed for that purpose. As before, $\n$ is random noise. In order to estimate $\A$ with some confidence, we must transmit as least as many pilot vectors as there are columns in $\A$. In addition we must assume that the realization of $\A$ does not change before all pilots have been transmitted. This assumption typically holds in flat, block-fading MIMO systems\cite{4548419, 1597555,5340650}.
With multiple transmitted pilots, model (\ref{linmod2}) can be written in matrix form as
\begin{equation}
\mathbf{Z=AS+N}.
\label{linmod3}
\end{equation}
If $\A$ is $m \times n$, then this model can be vectorized
into (Thm. 2, Ch. 2, \cite{Magn:Neud:1999})
\begin{equation}
\underbrace{\vec(\Z)}_{\y}=\underbrace{\left(\S^T \otimes \I_m \right)}_{\H}\underbrace{\vec(\A)}_{\x}+\underbrace{\vec(\N)}_{\n}.\label{linmod4}
\end{equation}
Here $\I_m$ denotes the $m \times m$ identity matrix,
the $\vec(\cdot)$ operator stacks the columns of a matrix into a column vector, and $\otimes$ denotes the Kronecker product.
Assuming that the channel ($\x$) and noise ($\n$) are independent and GM distributed, this is again a design problem as described in Section \ref{intro}, where the pilot matrix $\S$ is the design parameter.  A natural constraint, is to impose power limitations on the transmitted pilots, i.e. $\left\|\S\right\|^2_2\leq \gamma$. Together with the structure imposed by the Kronecker product, this then determines $\mathbb{H}$ in \eqref{theproblem}.

In \eqref{linmod4}, one may either assume that $\n$ is pure background noise, or that $\n$ represents noise \textit{and} interference. In the former case a Gaussian distribution may be justifiable, whereas in the latter a GM distribution may be more suitable \cite{1570066,276126}. As for the channel $\x$, a GM distribution can account for multiple fading situations. This can be useful, for example if the source is assumed to transmit from multiple locations. Then, the commonly used Rice distribution is unlikely to accurately capture the channel statistics associated with all transmit locations (especially so in urban areas). In fact, in \cite{Medboetal} it has been experimentally observed and reported that different transmit locations are indeed associated with different channel statistics. A GM distributed channel, with multiple modes, has the potential to capture this. 

The assumption that a channel realization can originate from several underlying distributions is not novel. For instance, all studies assuming channels governed by a Markov Model make this assumption, see e.g. \cite{350282,4607216} and the references therein. A GM is a special case of an Hidden Markov model, where subsequent observations are independent, rather than governed by a Markov process. In spite of this, to the best of our knowledge, pilot optimization for estimating channels governed by a GM distribution has not been considered in the literature.

\subsection {Gaussian Mixture distributions}
While aimed at minimizing the MSE, most optimization studies on linear precoders\cite{995062,1661076} or pilot signals\cite{4548419, 1597555,5340650} utilize only the first and second moments of the input distributions. Commonly, the underlying motivation is that a \textit{linear} MMSE (LMMSE) estimator is employed. The LMMSE estimator\footnote{Among all estimators which are linear (affine) in the observations, the LMMSE estimator obtains the smallest MSE.} only relies on first and second order statistics, which conveniently tends to simplify the associated matrix design problem. In fact, the desired matrix can often be obtained as the solution of a convex optimization problem. It is known, however, that the LMMSE estimator is optimal \textit{only} for the special case when the random signals $\x$ and $\n$ are both Gaussian. For all other cases, the LMMSE estimator is suboptimal. 

In practice, purely Gaussian inputs are rare. In general, the input distributions may be asymmetric, heavy tailed and/or multi modal.
A type of distribution that can accommodate all of these cases is the \textit{Gaussian Mixture} (GM) distribution.
In fact, a GM can in theory represent \textit{any} distribution with arbitrary accuracy \cite{Li99mixturedensity}, \cite{Sorenson1971465}. Therefore, in this work, we assume that the inputs are GM distributed as in \eqref{x_mix} and \eqref{n_mix}.
Notation \eqref{x_mix} should be read in the distributional sense, where $\mathbf{x}$ results from a composite experiment. First, source $k\in\mathcal{K}$ is activated with probability $p_k\geq0$, $\sum_{k\in\mathcal{K}} p_k =1$. Second, that source generates a Gaussian signal with distribution law $\mathcal{N}(\mathbf{u}^{(k)}_{\mathbf{x}},\mathbf{C}^{(k)}_{\mathbf{x}\mathbf{x}})$. For any realized $\mathbf{x}$, however, the underlying index $k$ is not observable. 
The noise $\mathbf{n}$ emerges in an entirely similar, but independent manner. $\mathcal{K}$ and $\mathcal{L}$ are index sets. In theory, it suffices that these sets are countable, but in practice they must be finite. Clearly when $\mathcal{K}$ and $\mathcal{L}$ are singletons, one falls back to the familiar case of Gaussian inputs. 


The mixture parameters, e.g. $(p_k,\mathbf{u}^{(k)}_{\mathbf{x}},\mathbf{C}^{(k)}_{\mathbf{x}\mathbf{x}})_{k \in \mathcal{K}}$, are rarely given a priori. Most often they must be estimated, which is generally a non-trivial task\cite{Li99mixturedensity,814639}. A common approach is to estimate the GM parameters from training data. The expectation maximization (EM) algorithm is well suited, and much used, for that purpose\cite{151045,Ghahramani93,mackay03}. Briefly, the algorithm relies on observations drawn from the distribution we wish to parametrize, and some initial estimate of the parameters. The observations are used to update the parameters, iteratively, until convergence to a local maximum of the likelihood function. Because the resulting GM parameters depend on the initial estimates, the algorithm can alternatively be started from multiple initial estimates. This produces multiple sets of GM parameters, and each set can be assigned probabilities based on the training data. Our starting point is that the distributions in (\ref{x_mix}) and \eqref{n_mix} have resulted from such, or similar model fitting.

Model \eqref{linmod} with GM inputs \eqref{x_mix} and \eqref{n_mix} is quite generic, and we have indicated two signal processing applications where the matrix design problem appears. The solution to that problem is, however, essentially application independent. It should therefore be of interest to a wide audience. To the best of our knowledge, it has not been pursued in the literature.

\section{An illustrative example}\label{A simple example}
\begin{figure}[htbp]
\begin{center}
\begin{pspicture}(8.5,4)
   \psline{->}(0,2)(4,2)
		\psline{->}(2,0)(2,4)
   \pscircle[ linestyle=dashed](1,2){0.5}
   \pscircle[ linestyle=dashed](1,2){0.2}
   \pscircle(1,2){0.6}
   \pscircle(1,2){0.4}
      \pscircle[ linestyle=dashed](3,2){0.5}
   \pscircle[ linestyle=dashed](3,2){0.2}
   \pscircle(3,2){0.6}
   \pscircle(3,2){0.4}
   \put(0.2,3.6){(a)}
      \psline{->}(4.5,2)(8.5,2)
		\psline{->}(6.5,0)(6.5,4)
   \pscircle[ linestyle=dashed](5.5,2){0.5}
   \pscircle[ linestyle=dashed](5.5,2){0.2}
   \pscircle(6.5,3){0.6}
   \pscircle(6.5,3){0.4}
      \pscircle[ linestyle=dashed](7.5,2){0.5}
   \pscircle[ linestyle=dashed](7.5,2){0.2}
   \pscircle(6.5,1){0.6}
   \pscircle(6.5,1){0.4}
   \put(4.7,3.6){(b)}
\end{pspicture}
\caption{(a): Densities without any rotation. (b): The effect of rotating $\mathbf{x}$ by $\pi/2$.}\label{fig:2}
\end{center}
\end{figure}
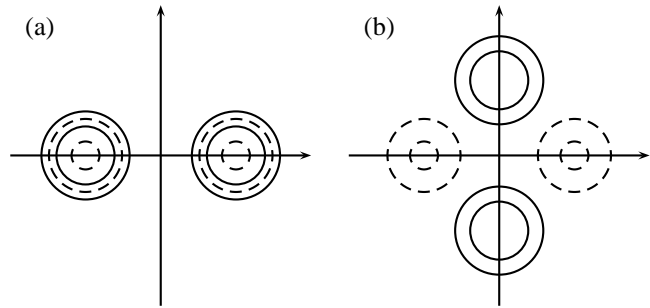

\begin{figure}
\includegraphics[scale=0.55]{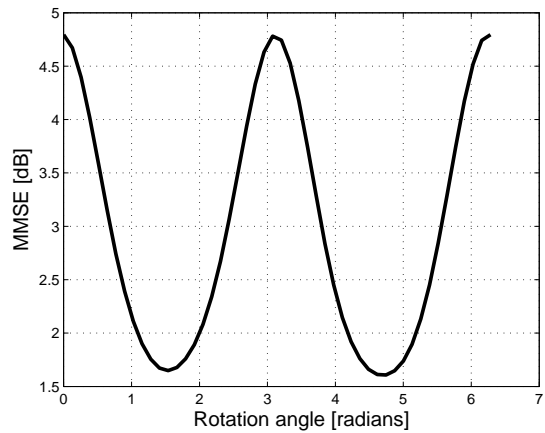}
\vspace{-0.3cm}
\caption{MMSE versus rotation angle.}
\label{fig:3}
\end{figure}

We start by studying a special instance of the matrix design problem, where the MMSE for \textit{all} $\H \in \mathbb{H}$ can be plotted. In general this is not possible, but the following simple example reveals some fundamental properties of the problem. 
Assume that we wish to design a precoder, as in \eqref{linmod5}, where $\B=\I_2$ and $\mathbf{F}$ is restricted to be an orthogonal matrix. 
Equation (\ref{linmod5}) then simplifies to $\mathbf{y}=\F\mathbf{x}+\mathbf{n}$, where $\F$ only \textit{rotates} $\x$.  Further, let $\mathbf{x}$ and $\mathbf{n}$ be independent and identically GM distributed as
\begin{align}\frac{1}{2}\mathcal{N}\left(\alpha\e_x,\mathbf{I}_2 \right)+\frac{1}{2}\mathcal{N}\left(-\alpha\e_x,\mathbf{I}_2\right),\nonumber
\end{align}
where $\alpha$ is a scalar and $\e_x$ is the unit vector along the $x$-axis.
Assume initially that $\F=\I_2$, which corresponds to no rotation. 
In this case, Figure \ref{fig:2}(a) illustrates the densities of $\mathbf{Fx}$ (full circles) and $\mathbf{n}$ (dashed circles), when seen from above. They are identical and sit on top of each other.
Now, with a precoder that rotates $\mathbf{x}$ by $\pi/2$, the densities of $\F\x$ and $\n$ will look like in Figure \ref{fig:2}(b). The latter configuration is preferable from an estimation viewpoint.
This is clear from figure \ref{fig:3}, where the MMSE is displayed as a function of all rotation angles between 0 and $2\pi$ (with $\alpha=2$). 
As can be seen, a significant gain can be obtained by rotating $\pi/2$ (or by $3\pi/2$). This gain is \textit{not} due to a particularly favorable signal-to-noise-ratio $\text{SNR}=E\left\|\F\x\right\|^2_2 /E\left\|\n\right\|^2_2 $; because $\mathbf{F}$ is orthogonal, the SNR remains equal for all rotation angles. The MMSE gain is instead due to a rotation producing a signal which \textit{tends} to be orthogonal to the noise. 

The above example is a special case of the matrix design problem, where $\H$ in (\ref{linmod}) is restricted to be orthogonal. It is clear that $\H$ plays a decisive role in how accurately $\x$ can be estimated. An almost equally important observation, however, is that the MMSE is not convex in $\H$. Hence, in general, we cannot expect that first order optimality (zero gradient) implies a \textit{global} minimum.
When studying the channel estimation problem further, we will see an implication of this non-convexity, which is perhaps not well known: In certain cases the MMSE of the channel estimate does not decrease with increasing pilot power. On the contrary, the MMSE may in fact increase.

In the next section we rewrite the original minimization problem into an equivalent but more compact maximization problem. Then, in Section, \ref{Kiefer-Wolfowitz precoding} we present a stochastic optimization approach which provides a solution.

\section{An equivalent maximization problem}\label{An equivalent optimization problem}
In order to propose a solution to the matrix design problem in (\ref{theproblem}), we first rewrite expression \eqref{defmmse}.
Using the results of \cite{flamchatterjeekansanen}, it follows that for model (\ref{linmod}), under independent GM inputs (\ref{x_mix}) and \eqref{n_mix}, and a fixed $\H$, the MMSE can be written as
\begin{align}
&E \left\{ \| \mathbf{x} -E\left\{\mathbf{x}|\mathbf{y}\right\} \|_{2}^{2}\right\} =\nonumber\\
&\hspace{0.2cm}\sum_{k}p_k \left( \text{tr}\left(\mathbf{C}^{(k)}_{\mathbf{x}\mathbf{x}}\right)+\left\|\mathbf{u}^{(k)}_{\mathbf{x}}\right\|^2_2 \right)
-\int \left\|\mathbf{u}_{\mathbf{x}|\mathbf{y}}\right\|^2_2 f(\mathbf{y})d\mathbf{y}. \label{mmse}
\end{align}
In (\ref{mmse}), $\text{tr}(\cdot)$ denotes the trace operator and $f(\y)$ is a (GM) probability density function
\begin{align}
f(\y)=\sum_{k,l}p_k q_l f^{(k,l)}(\y)\label{firsteq},
\end{align}
where
\begin{align}
f^{(k,l)}(\y)&=\frac{e^{-\frac{1}{2}\left(\y-\u^{(k,l)}_{\y}\right)^T \C^{-(k,l)}_{\y\y}\left(\y-\u^{(k,l)}_{\y}\right)}}{\left(2\pi\right)^{\frac{M}{2}}\left|\C^{(k,l)}_{\y\y}\right|^{\frac{1}{2}}},\label{dens}\\
\u^{(k,l)}_{\y}&=\H\u^{(k)}_{\x} +\u^{(l)}_{\n},\\
\C^{(k,l)}_{\y\y}&=\H\C^{(k)}_{\x\x} \H^T +\C^{(l)}_{\n\n}.
\end{align}
In \eqref{dens}, $(\cdot)^T$ denotes transposition, $\left|\cdot\right|$ denotes the determinant,  $\C^{-(k,l)}_{\y\y}$ is short for $(\C^{(k,l)}_{\y\y})^{-1}$ and $M$ is the length of $\y$.
The MMSE estimator $\mathbf{u}_{\mathbf{x}|\mathbf{y}}$ in (\ref{mmse}) can be written as
\begin{align}
\mathbf{u}_{\mathbf{x}|\mathbf{y}} = \frac{\sum_{k,l}p_k q_l f^{(k,l)}(\y)\mathbf{u}^{(k,l)}_{\mathbf{x}|\mathbf{y}}}{f(\y)}, \label{mmseest}
\end{align}
where
\begin{align}
\u^{(k,l)}_{\x|\y} &=  \u^{(k)}_{\x} + \C^{(k)}_{\x\x}\H^T \C^{-(k,l)}_{\y\y}\left(\y-\u^{(k,l)}_{\y}\right).
\label{comppostmean}
\end{align}
In what follows, it is convenient to define
 \begin{align}
G(\mathbf{H},\mathbf{y})\triangleq \left\|\mathbf{u}_{\mathbf{x}|\mathbf{y}}\right\|_2^2.\label{GHy}
\end{align}
This notation emphasizes that the squared norm of the MMSE estimate depends on both $\H$ and the observation $\y$.  In (\ref{mmse}), only the integral depends on $\H$.
Exploiting this, and using \eqref{GHy}, the minimizer of the MMSE is that $\mathbf{H}$ which maximizes
\begin{align}
\int G(\mathbf{H},\mathbf{y}) f(\mathbf{y})d\mathbf{y} =E\left[G(\mathbf{H},\mathbf{y})\right]=:g(\mathbf{H}), \label{objective}
\end{align}
subject to
\begin{align}
\hspace{0.2cm} \H\in \mathbb{H}. \label{constraint}
\end{align}
The integral in (\ref{objective}) cannot be evaluated analytically, even for a fixed and known $\H$ \cite{flamchatterjeekansanen}.
Moreover, as the example in Section \ref{A simple example} illustrated, the MMSE is generally not convex in $\H$, which implies that $g(\H)$ is generally not concave. Hence, any optimization method that merely aims at first order optimality, does in general not produce a global maximizer 
for $g(\H)$. Finally, as argued in the appendix, neither first or second order derivatives of $g(\H)$ w.r.t. $\H$ can be computed exactly, 
and accurate approximations of these are hard to obtain.

\section{The Robbins-Monro solution}\label{Kiefer-Wolfowitz precoding}
The above observations suggest that a sampling based approach is the only viable option. The problem of maximizing a non-analytical expectation $E\left[G(\mathbf{H},\mathbf{y})\right]$, over a parameter $\H$, falls under the umbrella of stochastic optimization. In particular, for our problem, the Robbins-Monro algorithm \cite{robbins1951, kushner2003stochastic}, can be used to move iteratively from a judicially chosen initial matrix $\mathbf{H}_0$ to a local maximizer $\mathbf{H}^*$. The philosophy is to update the current matrix $\H$ using the gradient of the MMSE. Since the gradient cannot be calculated, one instead relies on a stochastic approximation. Translated to our problem, the idea is briefly as follows. 
Although (\ref{objective}) cannot be computed analytically, it can be estimated from  independent sample vectors $\left\{\y_i = \H\x_i +\n_i\right\}^N_{i=1}$, as
\begin{align}
g(\H)\approx\frac{1}{N}\sum_i \left\|\mathbf{u}_{\mathbf{x}|\mathbf{y_i}}\right\|^2_2.\label{estm}
\end{align}
The derivative of (\ref{estm}) w.r.t. $\H$ represents an approximation of  $\frac{\partial g(\mathbf{H})}{\partial \mathbf{H}}$, which can be used to update the current $\H$. Each update is then projected onto the feasible set $\mathbb{H}$. 
This is the core idea of the much celebrated Robbins-Monro algorithm\cite{robbins1951}. 
In our context, the algorithm can be outlined as follows.
\begin{itemize}
	\item Let the current matrix be $\mathbf{H}_r$.
	\item Draw at random $(\x,\n)$ and compute $\y=\H_r \x + \n$.
		 \item Calculate the update direction as
	\begin{align}
\B_r= \frac{\partial G(\mathbf{H}_r,\mathbf{y})}{\partial \mathbf{H}_r},\label{Br}
\end{align}
and
	\begin{align}
\mathbf{W}_r=\mathbf{H}_r + \epsilon_r \B_r\label{candidate},
\end{align}
where $\left\{\epsilon_r\right\}^{\infty}_{r=1}$ is an infinite sequence of step sizes satisfying $\epsilon_r >0 $, $\epsilon_r \rightarrow 0$ and $\sum^{\infty}_{r=1} \epsilon_r = \infty$. 
\item Update the matrix as
\begin{align}
\mathbf{H}_{r+1}=\pi_{\mathbb{H}}\left(\mathbf{W}_r\right), \label{update}
\end{align}
where $\pi_{\mathbb{H}}(\cdot)$ represents the projection onto the set of permissible matrices $\mathbb{H}$. 
\item Repeat all steps until convergence.
\end{itemize}

\subsection{Remarks on the Robbins-Monro algorithm}
Recall that the input statistics \eqref{x_mix}, \eqref{n_mix} are assumed known. Therefore, in a design phase, it is reasonable to assume that the inputs $(\x,\n)$ can be sampled to \textit{compute} $\y$, as indicated in the second step of the algorithm. The alternative would be to sample $\y$ directly, leaving the underlying $\x$ and $\n$ unknown. The first approach is preferred because it guarantees that \eqref{Br} becomes an unbiased estimate of $\frac{\partial g(\H)}{\partial \H}$, whereas the alternative does not. This important point is fully explained in the appendix. In general, the Robbins-Monro procedure does not \textit{require} observing the input realizations $(\x,\n)$. The algorithm converges also when only outputs $\y$ are available. For this reason we write \eqref{Br} in terms of $\y$, but for our implementation we will assume that the underlying inputs $(\x,\n)$ are fully known.  

Because the gradient direction $\frac{\partial G(\mathbf{H},\mathbf{y})}{\partial \mathbf{H}}$ in \eqref{Br} is central in the algorithm, its closed form expression is derived in the appendix.  Observe that $\frac{\partial G(\mathbf{H},\mathbf{y})}{\partial \mathbf{H}}$ is random because it relies on a random realization of $\y$. Specifically it is a stochastic approximation of $\frac{\partial g(\H)}{\partial \H}$ based on a \textit{single} observation vector $\y$. Instantaneously, it may even point in directions opposite to $\frac{\partial g(\H)}{\partial \H}$. In order to increase the likelihood of a beneficial update, one can alternatively compute the gradient as an average based on multiple $\y$'s, as suggested in (\ref{estm}). Then
	\begin{align}
&\B_r=\frac{1}{N}\sum^N_{i=1} \frac{\partial G(\mathbf{H}_r,\mathbf{y}_i)}{\partial \mathbf{H}_r}.\nonumber
\end{align}
In our implementation of the algorithm, however, we do not do this. In fact, it was recognized by Robbins and Monro, that choosing $N$ large is generally inefficient. The reason is that $\H_r$ is only intermediate in the calculations, and as argued in the appendix, regardless of the value of  $N$, the update direction can be chosen such that it coincides with $\frac{\partial g(\mathbf{H})}{\partial \mathbf{H}}$ \textit{in expectation}. 

The Robbins-Monro procedure does not rely on the existence of $\frac{\partial G(\mathbf{H},\mathbf{y})}{\partial \mathbf{H}}$ at all points. If this derivative is discontinuous, one can instead use any of its sub-gradients; all of which are defined. Consequently, if the local maximum towards which the algorithm converges has a discontinuous derivative, then the algorithm will oscillate around this point. Due to the decaying step sizes, however, the oscillations will eventually become infinitesimal, and for all practical purposes, the system comes to rest.

Convergence towards a local optimum is guaranteed only as $r \rightarrow \infty$ \cite{robbins1951, kushner2003stochastic}. Therefore, in theory, the algorithm must run forever in order to converge. The engineering solution, which tends to work well in practice, is to terminate the algorithm when $\left\|\H_{r+1}-\H_{r}\right\|_{2} < \gamma$, where $\gamma$ is a chosen threshold, or simply after a predefined number of iterations. Still, the associated running time may be non-negligible, and therefore the Robins-Monro procedure is best suited when the input signals are stationary.

In general, for other problems than considered here, it may happen that the functional form of $G(\H,\y)$ is unknown, even when its output can be observed for any $\H$ and $\y$. In this case, $\frac{\partial G(\mathbf{H},\mathbf{y})}{\partial \H}$ cannot be computed. Instead one may replace it by a \textit{finite difference approximation} \cite{1952}. In some cases, this may also be preferable even when the derivative can be computed;  Especially so if computing $\frac{\partial G(\mathbf{H},\mathbf{y})}{\partial \H}$ requires much effort. 
When finite difference approximation are used, the procedure is known as the Kiefer-Wolfowitz algorithm \cite{1952, kushner2003stochastic}. If the derivative can be computed, however, the Kiefer-Wolfowitz algorithm is associated with more uncertainty (larger variance) than the Robbins-Monro procedure. For the interested reader, the present paper extends \cite{Flamisita}, which considers Kiefer-Wolfowitz precoding.

\section{Numerical results}\label{Numerical results}
In this section we will study two specific examples. One is on linear precoding, the other is on pilot design for channel estimation. In conformance with much of the literature, we will use the \textit{normalized} MSE (NMSE) as performance measure\footnote{Assuming $\x$ to be a zero mean signal, the NMSE is never larger than 1 (zero dB). The reason is that the MMSE estimator, $\u_{\mathbf{x}|\mathbf{y}}$, will coincide with the prior mean of $\x$ only when the SNR tends to zero. Hence, the prior mean is a worst case estimate of $\x$, and the NMSE describes the relative gain over the worst case estimate.}. This is defined as
\begin{align}
\text{NMSE}=\frac{E \left\{ \| \mathbf{x} -\u_{\mathbf{x}|\mathbf{y}} \|_{2}^{2}\right\}}{E \left\{ \| \mathbf{x} \|_{2}^{2}\right\}}. \nonumber
\end{align}
\subsection {Precoder design}
Here we study the performance of a Robbins-Monro precoder. 
As in the simple example of Section \ref{A simple example}, we restrict the precoder to be orthogonal. Thus, the norm of precoded signal is equal to that of the non-precoded signal. 
For the current example we choose the following parameters.
  \begin{itemize}
	\item $\mathbf{B}=\mathbf{I}_2$.
	\item $\mathbf{x}$ is GM distributed with parameters
		 \begin{align}
		 &p_k = 1/4,\text{ for } k=1...4\nonumber\\
 &\mathbf{u}^{(1)}_{\mathbf{x}}= \left [ \begin{array}{c} -10 \\ 10 \end{array} \right ] \text{ , }  \mathbf{u}^{(2)}_{\mathbf{x}}= \left [ \begin{array}{c} 10 \\ -10 \end{array} \right ],\nonumber\\
 &\mathbf{u}^{(3)}_{\mathbf{x}}= \left [ \begin{array}{c} 10 \\ 10 \end{array} \right ]\text{ , }  \mathbf{u}^{(4)}_{\mathbf{x}}= \left [ \begin{array}{c} -10 \\ -10 \end{array} \right ],\nonumber\\
  &\mathbf{C}^{(1)}_{\mathbf{x}\mathbf{x}}=\mathbf{C}^{(2)}_{\mathbf{x}\mathbf{x}}=\mathbf{C}^{(3)}_{\mathbf{x}\mathbf{x}}=\mathbf{C}^{(4)}_{\mathbf{x}\mathbf{x}}=\frac{1}{10}\mathbf{I}_2\nonumber.
 \end{align}
 	\item The noise is Gaussian and distributed as
\begin{align}
\mathbf{n}\sim\mathcal{N}\left(\left [ \begin{array}{c} 0 \\ 0\end{array} \right ],a\left [ \begin{array}{cc} 1 &  0\\ 0 & 0.1  \end{array} \right ]\right).
\end{align}
where $a$ is a scalar that can account for any chosen SNR=$ \text{tr} \left(\mathbf{C}_{\mathbf{x}\mathbf{x}}\right)/ \text{tr} \left(\mathbf{C}_{\mathbf{n}\mathbf{n}}\right)$, and $\mathbf{C}_{\mathbf{x}\mathbf{x}}$ and $\mathbf{C}_{\mathbf{n}\mathbf{n}}$ are the covariance matrices of $\x$ and $\n$ respectively.
\item We use $\mathbf{F}_0 = \mathbf{I_2}$ as the initial guess in the Robbins-Monro algorithm.
\item As stopping criterion we use: $\left\|\mathbf{F}_{r+1}-\mathbf{F}_{r}\right\|_2<10^{-4}$.
\end{itemize}

Because we have assumed $\B=\I_2$, we may in the Robbins-Monro algorithm of Section \ref {Kiefer-Wolfowitz precoding} simply replace all $\H_r$ with the precoder $\F_r$. For the projection in (\ref{update}), we choose the \textit{nearest orthogonal matrix}. This projection is the solution to the following optimization problem.
\begin{align}
\mathbf{F}_{r+1}&=\text{arg}\min_{\mathbf{F} \in \mathbb{O}}\left\|\mathbf{F}-\mathbf{W}_r \right\|^2_{2}\nonumber
\end{align}
where $\mathbb{O}$ is the set of orthogonal matrices. 
The solution is particularly simple, and exploits the singular value decomposition:
\begin{align}
&\mathbf{W}_r=\mathbf{U}\mathbf{D}\mathbf{V}^T \Rightarrow \mathbf{F}_{r+1} = \mathbf{U}\mathbf{V}^T.\nonumber
\end{align} 

\begin{figure}
\begin{center}
\includegraphics[scale=0.6]{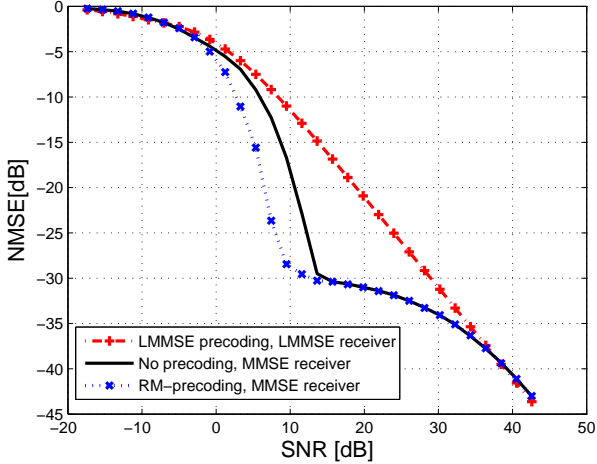}
\vspace{-0.2cm}
\caption{The NMSE with and without precoding.}
\label{fig:4}
\end{center}
\end{figure}

Figure \ref{fig:4} displays the NMSE with and without precoding, for increasing SNR levels. As can be seen, Robbins-Monro precoding provides a significant NMSE gain, especially at SNR levels between 0 and 10dB. Observe that the common approach of using the LMMSE estimator (and its corresponding precoder) is highly suboptimal at intermediate SNR levels. In fact, it is much worse than doing MMSE estimation without any precoding.

The above example indicates that our method generates a reasonable precoder, for \textit{particular} GM input distributions. Admittedly, there exists input statistics for which the gain is much less significant. However, in all simulations we have carried out, the clear tendency is that a Robbins-Monro precoder/ MMSE receiver outperforms the LMMSE precoder/LMMSE receiver at intermediate SNR levels. 

\subsection {Pilot design for channel estimation} \label{Pilot design for channel estimation}
Also for the channel estimation problem, the Robbins-Monro pilot matrix/ MMSE receiver outperforms the LMMSE pilot matrix/LMMSE receiver at intermediate SNR levels. In the next example we choose parameters in order to highlight this, and one additional property. That property is a direct consequence of the non-convex nature of the MMSE. We believe it to be of interest, but not well known.
The starting point is the channel estimation problem in (\ref{linmod3}), where we assume that all matrices are $2\times 2$. In the corresponding vectorized model
\begin{equation}
\underbrace{\vec(\Z)}_{\y}=\underbrace{\left(\S^T \otimes \I_2\right)}_{\H}\underbrace{\vec(\A)}_{\x}+\underbrace{\vec(\N)}_{\n},\label{vecvers}
\end{equation}
we assume the following parameters.
  \begin{itemize}
  \item The vectorized channel, $\x$, is distributed as $\mathcal{N}\left(\mathbf{0},\I_4\right)$.
	\item The vectorized noise, $\mathbf{n}$, is GM distributed with parameters
		 \begin{align}
		 &q_l = 1/2,\text{ for } l=1,2,\nonumber\\
 &\mathbf{u}^{(1)}_{\mathbf{n}}= -\mathbf{u}^{(2)}_{\mathbf{n}}=5\left [ \begin{array}{cccc} 1 & 1 & 1 & 1 \end{array} \right ]^T, \nonumber\\
  &\C^{(1)}_{\n\n}=\C^{(2)}_{\n\n}=\frac{1}{2}\I_4
 \end{align}
 	\item As constraint we impose that $\left\|\S\right\|^2_2=\alpha$,
where $\alpha$ is a positive scalar that can account for any chosen pilot power, and therefore also any SNR=$\left\|\S\right\|^2_2 /\text{tr} \left(\C_{\n\n}\right)$.
\item Stopping criterion: $\left\|\mathbf{S}_{r+1}-\mathbf{S}_{r}\right\|_2<10^{-4}$.
\end{itemize}
As starting point for the Robbins-Monro algorithm we set $\S$ equal to a scaled identity matrix satisfying the power constraint. During the iterations we rely on the following simple projection:
If the candidate pilot matrix has power $\left\|\S\right\|^2_2=\gamma$, then $\S\rightarrow \sqrt{\frac{\alpha}{\gamma}}\S$. Thus, if the pilot matrix does not use the entire power budget, the magnitude of all its elements are increased. Similarly, if pilot matrix has to large power, the magnitude of its elements are decreased. 

\begin{figure}
\begin{center}
\includegraphics[scale=0.6]{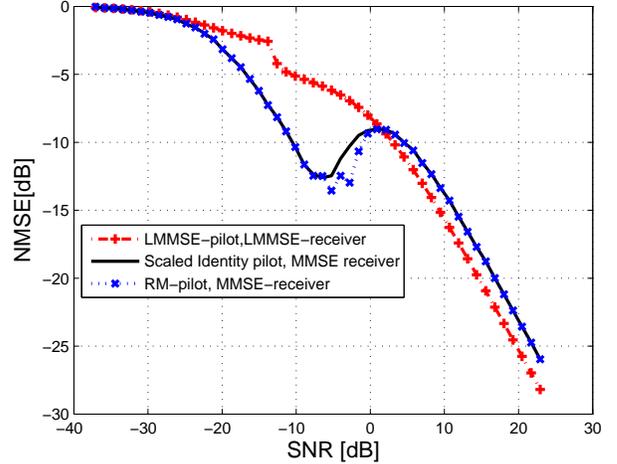}
\vspace{-0.5cm}
\caption{The NMSE as a function of pilot power.}
\label{fig:5}
\end{center}
\end{figure}

Figure \ref{fig:5} shows the estimation accuracy for increasing SNR (increasing values of $\alpha$). 
It can bee seen that our method outperforms the commonly used LMMSE channel estimator/LMMSE pilot matrix at intermediate SNRs. In fact, for the same range of SNRs, the latter is much than transmitting a scaled identity pilot matrix and using the MMSE estimator. As the SNR increases, however, it is known that the LMMSE estimator becomes optimal\cite{flamchatterjeekansanen}. The performance gap between our approach and the LMMSE estimator at high SNR indicates that a scaled identity pilot matrix is a local optimum that the Robbins-Monro algorithm does not easily escape from. The most striking observation in figure \ref{fig:5}, however, is that the channel estimates may become \textit{worse} by increasing the pilot power!
This is not an artifact of the Robbins-Monro algorithm; the same tendency is seen when a scaled identity (satisfying the same power constraint) is used as the pilot matrix.
 
 \subsection{Increased pilot power $\neq$ improved channel estimates}
\begin{figure}
\begin{center}
\includegraphics[scale=0.6]{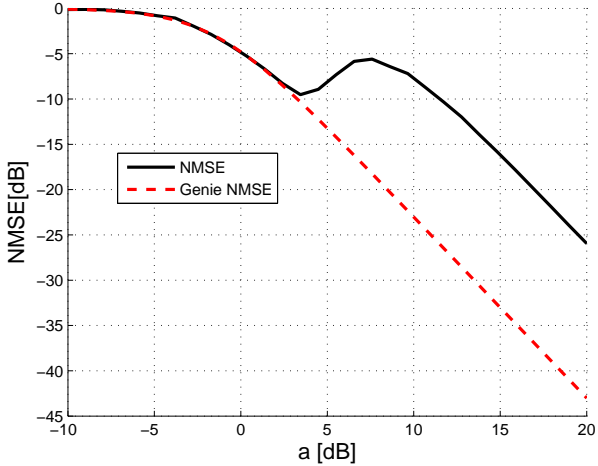}
\vspace{-0.5cm}
\caption{The NMSE as a function of the scalar $a$.}
\label{fig:6}
\end{center}
\end{figure}

We believe that the above phenomenon is not well known, and that it deserves to be explained. In order to visualize what happens, we will consider an example of smaller dimensions, but with similar properties as in the previous subsection. Specifically, we will assume that the unknown channel matrix $\A$ is $2 \times 1$, and that the pilot signal is just a scalar, $\s=a$. Then, using \eqref{linmod4} it follows that $\H=a\I_2$. Thus, in this setup, we do not \textit{optimize} anything, we only try to explain the NMSE as function of different values for $a$. We assume the following parameters:
  \begin{itemize}
  \item $\H=a\I_2$, where $a$ is a scalar that we can vary.
  \item The signal (channel) $\x$ is distributed as $\mathcal{N}\left(\mathbf{0},\I_2\right)$.
	\item The noise $\mathbf{n}$ is GM distributed with parameters
		 \begin{align}
		 &q_l = 1/2,\text{ for } l=1,2,\nonumber\\
 &\mathbf{u}^{(1)}_{\mathbf{n}}= -\mathbf{u}^{(2)}_{\mathbf{n}}=5\left [ \begin{array}{cccc} 1 & 1 \end{array} \right ]^T, \nonumber\\
  &\C^{(1)}_{\n\n}=\C^{(2)}_{\n\n}=\frac{1}{2}\I_2
 \end{align}
\end{itemize}
In figure \ref{fig:6}, the NMSE is plotted as a function of increasing values for the scalar $a$. We observe the same tendency as in figure \ref{fig:5}: for increasing values of $a$ (corresponding to increasing pilot power in figure \ref{fig:5}), the NMSE may increase. In figure \ref{fig:6}, we also plot the NMSE that would be obtained by a genie aided estimator \cite{flamchatterjeekansanen}. Briefly, the genie aided estimator knows from which underlying Gaussian source the noise $\n$ originates for each observation $\y$. Accordingly it can always produce the MMSE estimate corresponding to a purely Gaussian model. The genie aided estimator can of course not be implemented in practice, but because it is much better informed than the MMSE estimator, it provides a lower bound on the NMSE. Yet, from figure \ref{fig:6} we see that for $a<3.45$ dB, the MMSE estimator is able to pin-point the correct noise component.
\begin{figure}
\begin{center}
\includegraphics[scale=0.6]{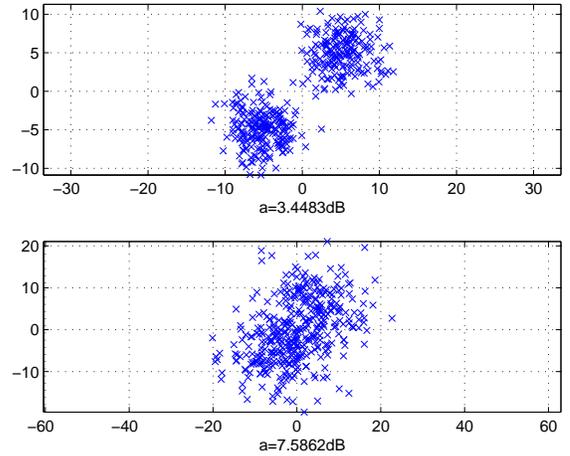}
\vspace{-0.5cm}
\caption{Sampled observations $\y$ for $a=3.45$ dB and $a=7.6$ dB.}
\label{fig:7}
\end{center}
\end{figure}
A plausible explanation is the following. 
For small $a$, almost all realizations of $\H\x=a\I_2 \x$ are close to the origin. Thus, observations $\y$ tend to appear in two distinct clusters; one cluster centered at each noise component. As a consequence, the active noise component can essentially always be identified. As $a$ grows, $\H\x=a\I_2 \x$ take values in an increasingly larger area, and the cluster borders approach each other. The value $a=3.45$ dB is the largest value for $a$ where the clusters can still be 'perfectly' separated. This value corresponds to the local minimum in figure \ref{fig:6}. Because we are considering 2-dimensional random vectors, we can actually visualize these clusters. The upper part of figure \ref{fig:7} shows how 400 independent $\y$'s form two nearby, but separable, clusters generated at $a=3.45$ dB. When $a$ grows beyond this level, the receiver faces a larger identification problem: it is harder to tell which noise component was active. The lower part of figure \ref{fig:7} shows 400 independent $\y$'s generated at $a=7.6$ dB. This value corresponds to the local maximum in figure \ref{fig:6}. Here the clusters largely overlap. As $a$ continues to grow, however, the average magnitude of a noise contribution becomes so small compared to the average magnitude of $a\I\x$, that near perfect recovery of $\x$ eventually becomes possible. 

Translated to the channel estimation problem in figure \ref{fig:5}, the interpretation is that there is a continuous range where increasing the pilot power is harmful. From  figure \ref{fig:5}, one observes that, unless one can spend an additional 15 dB (approximately) on pilot power, one will not improve from the local minimum solution. 

\section{Conclusion}
\label{concl}
We have provided a framework for solving the matrix design problem of the linear model under Gaussian mixture statistics. The study is motivated by two applications in signal processing. One
concerns the choice of error-reducing precoders; the other deals with
selection of pilot matrices for channel estimation. In either setting we use the Robbins-Monro procedure to arrive at a solution. Our
numerical results indicate improved estimation accuracy at intermediate SNR levels; markedly better
than those obtained by optimal design based on the LMMSE estimator.

Although the Robbins-Monro algorithm in theory only converges asymptotically, in practice we see that a hard stopping criterion may work well. The algorithm is still computationally demanding, and therefore best suited under stationary or near stationary settings. 

We have explored an interesting implication of the non-convexity of the MMSE; namely a case where spending more pilot power gives worse channel estimates. This phenomenon is not linked to the stochastic optimization procedure. It can be observed without optimizing $\H$ at all, and we have offered a plausible explanation. 

\section{Acknowledgments}
John T. Fl{\aa}m is supported by the Research Council of Norway under the
NORDITE/VERDIKT program, Project CROPS2 (Grant 181530/S10).

%
\appendix
This appendix derives a closed form expression for the gradient direction $\frac{\partial G\left(\H,\y\right)}{\partial \H}$ in \eqref{Br}, where $G\left(\H,\y\right)$ is defined through \eqref{mmseest}-\eqref{GHy}. To that end, it is worth observing that when optimizing $\H$ it is beneficial if the designer can draw samples directly from the inputs $\x$ and $\n$, and not only the output $\y$. In order to see why, 
assume in what follows that the order of derivation and integration can be interchanged such that
\begin{align}
 \frac{\partial g(\H)}{\partial \H}&=\frac{\partial}{\partial \H}\left(\int G(\mathbf{H},\mathbf{y}) f(\mathbf{y})d\mathbf{y}\right)\nonumber\\
 &= \int \frac{\partial}{\partial \H}\left[G\left(\mathbf{H},\mathbf{y}\right) f(\mathbf{y})d\mathbf{y}\right].\label{changeorder}
\end{align}
Now, if we can only observe outputs $\y$, we have
\begin{align}
&E\left(\frac{\partial G\left(\H,\y\right)}{\partial \H}\right)
= \int \frac{\partial G\left(\mathbf{H},\mathbf{y}\right)}{\partial \H} f(\mathbf{y})d\mathbf{y}\neq \frac{\partial g(\H)}{\partial \H}\nonumber.
\end{align}
Hence, in this case, the update direction $\frac{\partial G\left(\H,\y\right)}{\partial \H}$ is \textit{not} an unbiased estimator of the gradient $\frac{\partial g(\H)}{\partial \H}$. 
In contrast, assume that we can draw inputs $(\x,\n)$, and define 
\begin{align}
\left\|\mathbf{u}_{\mathbf{x}|\mathbf{\H\x+\n}}\right\|^2_2=\breve{G}\left(\H\x +\n\right),\nonumber
\end{align}
then
\begin{align}
&E\left( \frac{\partial \breve{G}\left(\H\x +\n\right)}{\partial \H}\right)=\iint \frac{\partial \breve{G}\left(\H\x +\n\right)}{\partial \H} f(\mathbf{x})d\mathbf{x}f(\mathbf{n})d\mathbf{n}\nonumber\\
&= \iint \frac{\partial}{\partial \H}\left[\breve{G}\left(\H\x +\n\right) f(\mathbf{x})d\mathbf{x}f(\mathbf{n})d\mathbf{n}\right]=\frac{\partial g(\H)}{\partial \H}\nonumber.
\end{align}
Here, the second equality holds because $f(\mathbf{x})d\mathbf{x}f(\mathbf{n})d\mathbf{n}$ is independent of $\H$. Hence, $\frac{\partial \breve{G}\left(\H\x +\n\right)}{\partial \H}$ \textit{is} an unbiased estimator of $\frac{\partial g(\H)}{\partial \H}$ , which is of course desirable.
Because it is beneficial to sample $\x$ and $\n$, rather than just $\y$, we will assume here that the designer can do this. In practice, this implies that the optimization of $\H$ is done off line, as preparation for the subsequent estimation.


In what follows, we will prove that interchanging the order of integration and derivation, as in \eqref{changeorder} is justifiable. We will derive a closed form expression for $\frac{\partial \breve{G}\left(\H\x +\n\right)}{\partial \H}$, and use this as the update direction in \eqref{Br}. Our strategy, however, will be to do this in the reverse order: First we compute the derivative, assuming that the change can be done, and then we show that that differentiation under the integral sign is justified. Although we assume knowledge of $(\x,\n)$ for each observed $\y$, we will write $\y$ instead of $\H\x +\n$, and $\frac{\partial G\left(\H,\y\right)}{\partial \H}$ instead of $\frac{\partial \breve{G}\left(\H\x +\n\right)}{\partial \H}$, simply to save space.

Using (\ref{mmseest}), $\frac{\partial G\left(\H,\y\right)}{\partial \H}$ can be written as
\begin{align}
 \sum_{k,l,r,s}p_k q_l p_r q_s\frac{\partial}{\partial \H}\left( \frac{  f^{(k,l)}(\y) f^{(r,s)}(\y){\u^{(k,l)}_{\x|\y}}^{T}\u^{(r,s)}_{\x|\y}}{\left(\sum_{k,l}p_k q_l f^{(k,l)}(\y)\right)^2}\right). \label{outerprod}
\end{align}
In order to compute \eqref{outerprod}, we make use of the following theorem \cite{Magn:Neud:1999}.
\begin{theorem}\label{th} For a scalar function, $\phi(\H)$, of a matrix argument, the \textit{differential} has the form
\begin{align}
 &d(\phi)= \text{tr}\left(\Q^T d(\H)\right)=\vec(\Q)^T \vec(d\H),\nonumber\\
&\text{where }\hspace{0.2cm} \Q=\frac{\partial \phi}{\partial \H}. \nonumber
  \end{align}
 \end{theorem} 
In our case, we take $\phi(\H)$ to be the expression in the large parenthesis of \eqref{outerprod}. We will identify its differential, and exploit Theorem \ref{th} in order to obtain the derivative. For that purpose, it is convenient to define
\begin{align}
f^{k,l,r,s}&=f^{(k,l)}(\y) f^{(r,s)}(\y),\\
z^{k,l,r,s}&={\u^{(k,l)}_{\x|\y}}^{T}\u^{(r,s)}_{\x|\y},\\
t&=\sum_{k,l}p_k q_l f^{(k,l)}(\y).
\end{align}
Using these, the derivative in (\ref{outerprod}), can then be compactly written as $\frac{\partial}{\partial \H}\left(\frac{f^{k,l,r,s}z^{k,l,r,s}}{t^2}\right)$. Using the chain rule, the differential of this fraction is
\begin{align}
&d(\phi)=d\left(\frac{f^{k,l,r,s}z^{k,l,r,s}}{t^2}\right)=- \frac{2f^{k,l,r,s}z^{k,l,r,s}d(t)}{t^3}\nonumber\\
&\hspace{1cm}+\frac{d(f^{k,l,r,s})z^{k,l,r,s}+d(z^{k,l,r,s})f^{k,l,r,s} }{t^2}.\label{simple}
\end{align}
Thus, we must identify the differentials $d(f^{k,l,r,s}), d(z^{k,l,r,s})$ and $d(t)$, which we do in next. \textbf{Notation:} we will in the remainder of this appendix use  $\left\langle \cdot\right\rangle$ to compactly denote the trace operator.\\
\\
\underline{Computing $d(f^{k,l,r,s})$}
\begin{align}
&d(f^{k,l,r,s})=d\left(f^{(k,l)}(\y) f^{(r,s)}(\y)\right)\nonumber\\
&=d\left(f^{(k,l)}(\y)\right) f^{(r,s)}(\y)+f^{(k,l)}(\y) d\left(f^{(r,s)}(\y)\right).\label{0}
\end{align}
The differential $d\left(f^{(k,l)}(\y)\right)$, is a differential of a Gaussian probability density function. In our case it depends on the indexes $(k,l)$, but in order to enhance readability, we will disregard these indexes in what follows. Thus, for now, we will use equations (\ref{firsteq})-(\ref{comppostmean}) with all indexes removed, and reincorporate the indexes when needed. In addition, we will for now disregard the constant factor $(2\pi)^{-\frac{M}{2}}$ in (\ref{dens}). Hence, instead of considering the differential $d\left(f^{(k,l)}(\y)\right)$, we therefore consider $d\left(\left|\C_{\y\y}\right|^{-\frac{1}{2}}g(\y)\right)$, where $g(\y)=e^{-\frac{1}{2}\left(\y-\u_{\y}\right)^T \C^{-1}_{\y\y}\left(\y-\u_{\y}\right)}$. This can be written as
\begin{align}
&d\left(\left|\C_{\y\y}\right|^{-\frac{1}{2}}g(\y)\right)\nonumber\\
&=d\left(\left|\C_{\y\y}\right|^{-\frac{1}{2}}\right)g(\y)+\left|\C_{\y\y}\right|^{-\frac{1}{2}}d\left(g(\y)\right)\nonumber\\
&=-\frac{g(\y)}{2}\left|\C_{\y\y}\right|^{-\frac{3}{2}}d\left(\left|\C_{\y\y}\right|\right)\nonumber\\
&\hspace{0.4cm}-\frac{g(\y)}{2}\left|\C_{\y\y}\right|^{-\frac{1}{2}}d\left(\left(\y-\u_{\y}\right)^T \C^{-1}_{\y\y}\left(\y-\u_{\y}\right)\right)\label{firstdiff}
\end{align}
In the second equality we have used the chain rule, and exploited that $g(\y)$ is an exponential function.
The first differential in (\ref{firstdiff}), provided $\C_{\y\y}$ is full rank, is (Theorem 1, ch. 8, of \cite{Magn:Neud:1999})
\begin{align}
 d\left(\left|\C_{\y\y}\right|\right)&=\left|\C_{\y\y}\right|\left\langle\C^{-1}_{\y\y}d\left(\C_{\y\y}\right)\right\rangle\label{1.1}\\
&=\left|\C_{\y\y}\right|\left\langle\C^{-1}_{\y\y}d\left(\H\C_{\x\x} \H^T +\C_{\n\n}\right)\right\rangle\nonumber\\
&=\left|\C_{\y\y}\right|\left\langle\C^{-1}_{\y\y}\left( d(\H)\C_{\x\x} \H^T +\H\C_{\x\x}d\left(\H^T\right)\right)\right\rangle\label{1.3}\\
&=\left|\C_{\y\y}\right|\left\langle\C_{\x\x} \H^T \C^{-1}_{\y\y} d(\H) + d(\H)\C_{\x\x} \H^T\C^{-1}_{\y\y}\right\rangle\label{1.4}\\
&=\left|\C_{\y\y}\right|\left\langle\C_{\x\x} \H^T \C^{-1}_{\y\y} d(\H) +\C_{\x\x} \H^T\C^{-1}_{\y\y} d(\H)\right\rangle\label{1.5}\\
&=2\left|\C_{\y\y}\right|\left\langle\C_{\x\x} \H^T \C^{-1}_{\y\y} d(\H)\right\rangle\nonumber\\
&=2\left|\C_{\y\y}\right|\left\langle\left(\C^{-1}_{\y\y}\H\C_{\x\x}\right)^T d(\H)\right\rangle.\label{1.8}
\end{align}
In (\ref{1.4}), we have rotated the first trace (done a cyclic permutation of the matrix product), and transposed the second trace. Because $\C_{\x\x}$ and $\C^{-1}_{\y\y}$ are symmetric, they are not affected by transposition. Moreover, $d(\H^T)= (d(\H))^T$. The trace operator is invariant to such rotations and transposition, and therefore these operations are justified. In (\ref{1.5}) we have rotated the second term. Such rotations and transpositions will be frequently employed throughout. Introducing $\w=\y-\u_{\y}$, the second differential of (\ref{firstdiff}) can be written
\begin{align}
&d\left(\w^T \C^{-1}_{\y\y} \w\right)=d\left( \left\langle\w^T \C^{-1}_{\y\y} \w\right\rangle\right)\nonumber\\
&=  \left\langle\w^T d\left(\C^{-1}_{\y\y}\right) \w\right\rangle+2 \left\langle\w^T \C^{-1}_{\y\y} d\left(\w\right)\right\rangle.\label{expterm}
\end{align}
The first term of (\ref{expterm}) is
\begin{align}
&\left\langle\w^T d\left(\C^{-1}_{\y\y}\right)\w\right\rangle\label{2.1}\\
&= -\left\langle\w^T \C^{-1}_{\y\y}d\left(\C_{\y\y}\right)\C^{-1}_{\y\y}\w \right\rangle \label{2.2}\\
&= -\left\langle\w^T \C^{-1}_{\y\y}\left(d(\H)\C_{\x\x} \H^T +\H\C_{\x\x}d\left(\H^T\right)\right)\C^{-1}_{\y\y}\w\right\rangle\nonumber\\
&= -\left\langle\C^{-1}_{\y\y}\w\w^T \C^{-1}_{\y\y}\left( d(\H)\C_{\x\x} \H^T +\H\C_{\x\x}d\left(\H^T\right)  \right)\right\rangle.\label{2.5}
\end{align}
Equation (\ref{2.2}) results from Theorem 3, ch. 8, of \cite{Magn:Neud:1999}. Observe that $\C^{-1}_{\y\y}\w\w^T \C^{-1}_{\y\y}$ in (\ref{2.5}) is a symmetric matrix, playing the same role as $\C^{-1}_{\y\y}$ in (\ref{1.3}). Therefore, we can utilize (\ref{1.8}) and conclude that
\begin{align}
\left\langle\w^T d\left(\C^{-1}_{\y\y}\right)\w\right\rangle=-2\left\langle\left( \C^{-1}_{\y\y}\w\w^T \C^{-1}_{\y\y}\H\C_{\x\x} \right)^T d(\H)\right\rangle.\label{2.6}
\end{align}
Recall that $\w=\H(\x-\u_{\x})+\n-\u_{\n}$.
The second term of (\ref{expterm}) can therefore be written as
\begin{align}
&2 \left\langle\w^T \C^{-1}_{\y\y} d\left(\w\right)\right\rangle\nonumber\\
& = 2 \left\langle\w^T \C^{-1}_{\y\y}  d\left(\H\right)\left(\x-\u_{\x} \right)\right\rangle\nonumber\\
& = 2 \left\langle\left(\x-\u_{\x} \right)\w^T \C^{-1}_{\y\y}  d\left(\H\right)\right\rangle\nonumber\\
& = 2 \left\langle\left( \C^{-1}_{\y\y} \w \left(\x-\u_{\x} \right)^T \right)^T d\left(\H\right)\right\rangle.\label{3.1}
\end{align}
Using (\ref{1.8}),(\ref{2.6}) and (\ref{3.1}), and inserting into (\ref{firstdiff}), we find that 
\begin{align}
&d\left(\left|\C_{\y\y}\right|^{-\frac{1}{2}}g(\y)\right)\nonumber\\
&=-g(\y)\left|\C_{\y\y}\right|^{-\frac{1}{2}}\left\langle\left(\C^{-1}_{\y\y}\H\C_{\x\x}\right)^T d(\H)\right\rangle\nonumber\\
&\hspace{0.4cm}+g(\y)\left|\C_{\y\y}\right|^{-\frac{1}{2}} \left\langle\left( \C^{-1}_{\y\y}\w\w^T \C^{-1}_{\y\y}\H\C_{\x\x} \right)^T d(\H)\right\rangle\nonumber\\
&\hspace{0.4cm}-g(\y)\left|\C_{\y\y}\right|^{-\frac{1}{2}} \left\langle\left(\C^{-1}_{\y\y} \w\left(\x-\u_{\x} \right)^T\right)^T d\left(\H\right)\right\rangle.\label{6}
\end{align}
If we define 
\begin{align}
\R^{(k,l)}&= \C^{-(k,l)}_{\y\y} \w^{(k,l)}\left(\x-\u^{(k)}_{\x} \right)^T \nonumber\\
&\hspace{0.4cm}+\C^{-(k,l)}_{\y\y}\left( \I -\w^{(k,l)} {\w^{(k,l)}}^T \C^{-(k,l)}_{\y\y}\right) \H\C^{(k)}_{\x\x},\nonumber
\end{align}
where $\w^{(k,l)}=\y-\u^{(k,l)}_{\y}$, and reincorporate the constant factor $(2\pi)^{-\frac{M}{2}}$, we now find that
\begin{align}
d\left(f^{(k,l)}(\y)\right)=-f^{(k,l)}(\y) \left\langle\left(\R^{(k,l)}\right)^T d\left(\H\right)\right\rangle.\label{fdiff}
\end{align}
Accordingly, (\ref{0}) becomes
\begin{align}
d\left(f^{k,l,r,s}\right)-\left\langle f^{k,l,r,s}\left( \R^{(k,l)} +\R^{(r,s)} \right)^T d(\H)\right\rangle.\label{firstdiff2}
\end{align}\\
\\
\underline{Computing $d(z^{k,l,r,s})$}
\begin{align}
d(z^{k,l,r,s})=&d\left({\u^{(k,l)}_{\x|\y}}^{T}\u^{(r,s)}_{\x|\y}\right)\nonumber\\
&=\left\langle d\left({\u^{(k,l)}_{\x|\y}}^{T}\u^{(r,s)}_{\x|\y}\right)\right\rangle\nonumber\\
&=\left\langle{\u^{(r,s)}_{\x|\y}}^Td\left({\u^{(k,l)}_{\x|\y}}\right)+{\u^{(k,l)}_{\x|\y}}^{T}d\left(\u^{(r,s)}_{\x|\y}\right)\right\rangle.\label{seconddiff}
\end{align}
Apart from a rearrangement of the indexes, equation (\ref{seconddiff}) contains two similar terms. Hence it suffices to compute one of them. Recalling that $\u^{(k,l)}_{\x|\y}$ is defined by \eqref{comppostmean}, we focus on the differential
\begin{align}
&\left\langle{\u^{(r,s)}_{\x|\y}}^Td\left({\u^{(k,l)}_{\x|\y}}\right)\right\rangle\nonumber\\
&=\left\langle{\u^{(r,s)}_{\x|\y}}^Td\left(\u^{(k)}_{\x} + \C^{(k)}_{\x\x}\H^T \C^{-(k,l)}_{\y\y}\w^{(k,l)}\right)\right\rangle\nonumber\\
&=\left\langle{\u^{(r,s)}_{\x|\y}}^T \C^{(k)}_{\x\x}d\left(\H^T\right) \C^{-(k,l)}_{\y\y}\w^{(k,l)}\right\rangle\label{z1}\\
&\hspace{0.2cm}+\left\langle{\u^{(r,s)}_{\x|\y}}^T \C^{(k)}_{\x\x}\H^T d\left(\C^{-(k,l)}_{\y\y}\right)\w^{(k,l)}\right\rangle\label{z2}\\
&\hspace{0.2cm}+\left\langle{\u^{(r,s)}_{\x|\y}}^T \C^{(k)}_{\x\x}\H^T \C^{-(k,l)}_{\y\y}d\left(\w^{(k,l)}\right)\right\rangle.\label{z3}
\end{align}
We will resolve this term by term. The first term, (\ref{z1}), reads
\begin{align}
&\left\langle{\u^{(r,s)}_{\x|\y}}^T \C^{(k)}_{\x\x}d\left(\H^T\right) \C^{-(k,l)}_{\y\y}\w^{(k,l)}\right\rangle\nonumber\\
&=\left\langle\left( \C^{-(k,l)}_{\y\y}\w^{(k,l)}{\u^{(r,s)}_{\x|\y}}^T\C^{(k)}_{\x\x}\right)^T d\left(\H\right) \right\rangle.\label{z11}
\end{align}
The second term, (\ref{z2}), can be written as
\begin{align}
&\left\langle{\u^{(r,s)}_{\x|\y}}^T \C^{(k)}_{\x\x}\H^T d\left(\C^{-(k,l)}_{\y\y}\right)\w^{(k,l)}\right\rangle\nonumber\\
&=-\left\langle{\u^{(r,s)}_{\x|\y}}^T \C^{(k)}_{\x\x}\H^T \C^{-(k,l)}_{\y\y}d\left(\C^{(k,l)}_{\y\y}\right)\C^{-(k,l)}_{\y\y}\w^{(k,l)}\right\rangle\nonumber\\
&=-\left\langle\underbrace{\C^{-(k,l)}_{\y\y}\w^{(k,l)}{\u^{(r,s)}_{\x|\y}}^T \C^{(k)}_{\x\x}\H^T \C^{-(k,l)}_{\y\y}}_{\C^{(k,l,r,s)}}d\left(\C^{(k,l)}_{\y\y}\right)\right\rangle\nonumber\\
&=-\left\langle\C^{(k,l,r,s)}\left( d\left(\H\right)\C^{(k)}_{\x\x} \H^T +\H\C^{(k)}_{\x\x}d\left(\H^T\right) \right) \right\rangle\nonumber\\
&=-\left\langle\C^{(k)}_{\x\x} \H^T \left(\C^{(k,l,r,s)} + {\C^{(k,l,r,s)}}^T\right) d(\H)\right\rangle\nonumber\\
&=-\left\langle\left(\left(\C^{(k,l,r,s)} + {\C^{(k,l,r,s)}}^T\right)\H\C^{(k)}_{\x\x}\right)^T d(\H)\right\rangle.\label{z22}
\end{align}
The third term, (\ref{z3}), reads
\begin{align}
&\left\langle{\u^{(r,s)}_{\x|\y}}^T \C^{(k)}_{\x\x}\H^T \C^{-(k,l)}_{\y\y}d\left(\w^{(k,l)}\right)\right\rangle\nonumber\\
&=\left\langle{\u^{(r,s)}_{\x|\y}}^T \C^{(k)}_{\x\x}\H^T \C^{-(k,l)}_{\y\y} d\left(\H\right)\left(\x-\u^{(k)}_{\x} \right)\right\rangle\nonumber\\
&=\left\langle\left(\C^{-(k,l)}_{\y\y}\H\C^{(k)}_{\x\x}{\u^{(r,s)}_{\x|\y}}\left(\x-\u^{(k)}_{\x} \right)^T \right)^T d\left(\H\right)\right\rangle.\label{z33}
\end{align}
Using (\ref{z11}),(\ref{z22}) and (\ref{z33}) we now define
\begin{align}
\D^{(k,l,r,s)}&= \C^{-(k,l)}_{\y\y}\w^{(k,l)}{\u^{(r,s)}_{\x|\y}}^T\C^{(k)}_{\x\x}\nonumber\\
&\hspace{0.4cm}-\left(\C^{(k,l,r,s)} + {\C^{(k,l,r,s)}}^T\right)\H\C^{(k)}_{\x\x} \nonumber\\
&\hspace{0.4cm}+ \C^{-(k,l)}_{\y\y}\H\C^{(k)}_{\x\x}{\u^{(r,s)}_{\x|\y}}\left(\x-\u^{(k)}_{\x} \right)^T.\nonumber
\end{align}
Due to its two similar terms, the differential in (\ref{seconddiff}) can then be written
\begin{align}
d(z^{k,l,r,s})&=\left\langle\left(\D^{(k,l,r,s)} + {\D^{(r,s,k,l)}}\right)^T d(\H)\right\rangle.\label{seconddiff2}
\end{align}\\
\\
\underline{Computing $d(t)$}
\begin{align}
d(t)&=d\left(\sum_{k,l}p_k q_l f^{(k,l)}(\y)\right)=\sum_{k,l}p_k q_l d\left(f^{(k,l)}(\y)\right)\nonumber\\
&=-\sum_{k,l}p_k q_l f^{(k,l)}(\y)\left\langle\left(\R^{(k,l)}\right)^T d(\H)\right\rangle.\label{thirddiff}
\end{align}
The last equation results immediately by employing (\ref{fdiff}).

\subsection{Computing the derivative}
Utilizing (\ref{firstdiff2}), (\ref{seconddiff2}) and (\ref{thirddiff}), the complete differential in (\ref{simple}) can now be written as
\begin{align}
&d(\phi)=d\left(\frac{f^{k,l,r,s}z^{k,l,r,s}}{t^2}\right)\nonumber\\
&=-\frac{\left\langle f^{k,l,r,s}\left( \R^{(k,l)} +\R^{(r,s)} \right)^T d(\H)\right\rangle z^{k,l,r,s}}{t^2}\nonumber\\
&+\frac{\left\langle\left(\D^{(k,l,r,s)} + {\D^{(r,s,k,l)}}\right)^T d(\H)\right\rangle f^{k,l,r,s}}{t^2}\nonumber\\
&+\frac{2 f^{k,l,r,s}z^{k,l,r,s}\sum_{k,l}p_k q_l f^{(k,l)}(\y)\left\langle\left(\R^{(k,l)}\right)^T d(\H)\right\rangle}{t^3}.\label{finaldiff}
\end{align}
In case of a precoder design problem, one makes the following substitutions: $\H=\B\F$ and $d(\H)=\B d(\F)$ throughout. 
In case of the pilot design problem \eqref{linmod4}, $\H$ must be substituted by $\S^T \otimes \I_m$. In addition, assuming that $\S$ is $n \times r$, one makes use of the fact that
\begin{align}
 \vec(d \H)&= \vec\left(d(\S^T) \otimes \I_m\right)\nonumber\\
& = \left(\I_n \otimes \K_{mr} \otimes \I_m \right)\left(\I_{rn} \otimes \vec(\I_m)\right)d(\vec(\S^T)).\nonumber 
\end{align}
Here $\K_{mr}$ is the Magnus and Neudecker commutation matrix \cite{Magn:Neud:1999}.
Theorem \ref{th} can then be easily applied to (\ref{finaldiff}), and identifying the derivative in (\ref{outerprod}) is therefore now straightforward. 

Finally, assume that $\H$ in (\ref{finaldiff}) is \textit{not} a function of some other matrix. Complactly defining $p_k q_l p_r q_s f^{k,l,r,s}=h^{k,l,r,s}$, and observing that $\sum_{k,l,r,s}h^{k,l,r,s}=t^2$, we find from equations \eqref{outerprod}, \eqref{finaldiff}, and Theorem \ref{th} that
\begin{align}
&\frac{\partial\breve{ G}\left(\H\x +\n\right)}{\partial \H}=\nonumber\\
&-\frac{\sum_{k,l,r,s} {h^{k,l,r,s}\left( \R^{(k,l)} +\R^{(r,s)} \right) z^{k,l,r,s}}}{\sum_{k,l,r,s} h^{k,l,r,s}}\nonumber\\
&+\frac{\sum_{k,l,r,s}h^{k,l,r,s}\left(\D^{(k,l,r,s)} + {\D^{(r,s,k,l)}}\right) }{\sum_{k,l,r,s}h^{k,l,r,s}}\nonumber\\
&+\frac{2 \sum_{k,l,r,s}h^{k,l,r,s}z^{k,l,r,s}\sum_{k,l}p_k q_l f^{(k,l)}(\y)\R^{(k,l)}}{\sum_{i,j}p_i q_j f^{(i,j)}(\y)\sum_{k,l,r,s}h^{k,l,r,s}}.\label{gaygaz}
\end{align}

\subsection{Interchanging the order of derivation and integration}
Recall, that interchanging the order of derivation and integration, as in (\ref{changeorder}), was until now only assumed valid. It derives from Lebesgue's Dominated Convergence Theorem that such a change \textit{is} valid if there exists a \textit{dominating} function $v(\cdot)$ satisfying
\begin{align}
 \left\|\frac{\partial \breve{G}\left(\H\x +\n\right)}{\partial \H}\right\|_2 \leq \left\|v(\H,\x,\n)\right\|_2\label{boundnorm}
 \end{align}
  and
  \begin{align}
\iint \left\|v(\H,\x,\n)\right\|_2 f(\x)f(\n)d\x d\n <\infty. \label{boundint}
\end{align}
Now consider \eqref{gaygaz} and define the function 
\begin{align}
w_{k,l,r,s}(\H,\x,\n)=-&\left( \R^{(k,l)} +\R^{(r,s)} \right) z^{k,l,r,s}\nonumber\\
&+\left(\D^{(k,l,r,s)} + {\D^{(r,s,k,l)}}\right) \nonumber\\
&+\frac{2 z^{k,l,r,s}\sum_{k,l}p_k q_l f^{(k,l)}(\y)\R^{(k,l)}}{\sum_{i,j}p_i q_j f^{(i,j)}(\y)}.\nonumber
\end{align}
Observe that
\begin{align}
\frac{\sum_{k,l,r,s}h^{k,l,r,s}w_{k,l,r,s}(\H,\x,\n)}{\sum_{k,l,r,s}h^{k,l,r,s}}=\frac{\partial\breve{ G}\left(\H\x +\n\right)}{\partial \H}\nonumber.
\end{align} 
Hence $\frac{\partial\breve{ G}\left(\H\x +\n\right)}{\partial \H}$ is a convex combination of the $w_{k,l,r,s}(\H,\x,\n)$'s, and therefore the function
\begin{align}
v(\H,\x,\n)=\sum_{k,l,r,s}\left\|w_{k,l,r,s}(\H,\x,\n)\right\|_2\nonumber
\end{align} 
 clearly satisfies (\ref{boundnorm}). 
We do not explcitly prove it her, but it can be verified that the integral
\begin{align}
\iint \sum_{k,l,r,s} \left\|w_{k,l,r,s}(\H,\x,\n)\right\|_2 f(\x)f(\n)d\x d\n\label{ssint} 
\end{align}
is bounded. 
Hence a dominating function exists, and the change of integration and derivation is justified. 
\subsection{First and second order derivatives of the objective function}
When trying to compute
\begin{align}
\iint\frac{\partial \breve{G}\left(\H\x +\n\right)}{\partial \H}f(\x)f(\n)d\x d\n \label{last},
\end{align}
the mixture densities in the denominators of (\ref{gaygaz}) will not simplify by substitutions. An entirely similar argument provides the reason for why \eqref{objective} cannot be computed analytically in the first place \cite{flamchatterjeekansanen}. Hence, (\ref{last}) cannot be computed analytically, and a closed form derivative of \eqref{objective} w.r.t $\H$ does not exist. Although not demonstrated here, a similar argument will hold also for the second order derivative. 

\bibliographystyle{IEEEbib}
\bibliography{strings}
\end{document}